\documentclass[11pt]{article}

\usepackage{latexsym,amsmath,amsthm,amssymb,graphicx,epsfig,multicol,graphicx,epstopdf}
\def\fullpage {
\addtolength{\topmargin}{-2 cm}
\addtolength{\oddsidemargin}{-1.4 cm} \addtolength{\textwidth}{+2.8 cm}
\addtolength{\textheight}{+3.3 cm}}
\fullpage
\newcommand{\ex}{\mbox{ex}}

\makeatletter
\def\thm@space@setup{\thm@preskip=10pt
\thm@postskip=5pt}
\makeatother
\newtheoremstyle{newstyle}
{} %Aboveskip
{} %Below skip
{\itshape} %Body font e.g.\mdseries,\bfseries,\scshape,\itshape
{} %Indent
{\bfseries} %Head font e.g.\bfseries,\scshape,\itshape
{.} %Punctuation after theorem header
{ } %Space after theorem header
{} %Heading

\theoremstyle{newstyle}
\newtheorem{theorem}{Theorem}[section]
\newtheorem{proposition}{Proposition}
\newtheorem{lemma}{Lemma}

\newtheorem{conj}{Conjecture}

\makeatletter
\newenvironment{pf}[1][\proofname]{\par
  \pushQED{\qed}%
  \normalfont \topsep0\p@\relax
  \trivlist
  \item[\hskip\labelsep\itshape
  #1\@addpunct{.}]\ignorespaces
}{%
  \popQED\endtrivlist\@endpefalse
}
\makeatother

\newcommand{\fp}{\mathbb{F}_p}

\begin{document}
 \date{}

\title{\vspace{-1cm}A counterexample to sparse removal}

\author{Craig Timmons \\
\small Department of Mathematics\\ [-0.8ex]
\small  University of California at San Diego\\ [-0.8ex]
\small California, U.S.A.\\
\small{\tt ctimmons@ucsd.edu} \\
\and
Jacques Verstraete \thanks{Research supported by NSF Grant DMS 1101489.}\\
\small Department of Mathematics\\[-0.8ex]
\small University of California at San Diego\\[-0.8ex]
\small California, U.S.A.\\
\small {\tt jverstraete@math.ucsd.edu}
}

\maketitle

\begin{abstract}
The {\em Tur\'{a}n number} of a graph $H$, denoted $\mbox{ex}(n,H)$, is the maximum number of edges
in an $n$-vertex graph with no subgraph isomorphic to $H$. Solymosi~\cite{Sol} conjectured that if $H$ is any graph and $\mbox{ex}(n,H) = O(n^{\alpha})$ where $\alpha > 1$, then any
$n$-vertex graph with the property that each edge lies in exactly one copy of $H$ has $o(n^{\alpha})$ edges. This can be viewed as conjecturing a possible extension of the removal lemma to sparse graphs, and is well-known to be true when $H$ is a non-bipartite graph, in particular when $H$ is a triangle, due to Ruzsa and Szemer\'{e}di~\cite{RS}.
Using Sidon sets we exhibit infinitely many bipartite graphs $H$ for which the conjecture is false.
\end{abstract}

%%%%%%%%%%%%%%%%%%%%%%%%%%%%%%%%%%%%%%%%%%%%%%%%%

\section{Introduction}\label{intro}

The Removal Lemma~\cite{EFR,CF} states that if $(G_n)_{n \in \mathbb N}$ is a sequence of graphs where $G_n$ has $n$ vertices, and $H$ is a $k$-vertex graph such that $G_n$ contains $o(n^k)$ subgraphs isomorphic to $H$,
then $o(n^2)$ edges may be deleted from $G_n$ to obtain an $H$-free graph. Ruzsa and Szemer\'{e}di~\cite{RS} established this result in the case that $H$ is a triangle as a consequence of
Szemer\'{e}di's Regularity Lemma~\cite{Sz}, and then applied the removal lemma to give a weak form of Roth's Theorem~\cite{Roth} on three-term arithmetic progressions. The removal lemma is a central tool in extremal combinatorics with many
applications~\cite{CF}.

\medskip

The removal lemma is not effective if $G_n$ has $o(n^2)$ edges. In particular, if $H$ is a $k$-vertex bipartite graph and $G_n$ contains $o(n^k)$ copies of $H$,
then it is a consequence of extremal graph theory that $G_n$ has $o(n^2)$ edges. The {\em Tur\'{a}n number} $\ex(n,H)$ is the maximum number of edges in an $n$-vertex
graph not containing $H$. When $H$ is a bipartite graph, the K\"{o}vari-S\'{o}s-Tur\'{a}n Theorem~\cite{KST} immediately shows $\ex(n,H) = O(n^{2 - 1/k})$, and more
precise results are available~\cite{AKS,FoxSudakov} dependent on the finer structure of $H$. Solymosi~\cite{Sol} conjectured an extension of the removal lemma in this regime
as a function of the Tur\'{a}n number of $H$ as follows. The {\em exponent} of a graph $H$, when it exists, is a real number $\alpha$ such that $\mbox{ex}(n,H) = \Theta(n^{\alpha})$ as $n \rightarrow \infty$.
Erd\H{o}s and Simonovits~\cite{ES} conjectured that every graph has an exponent, but this conjecture remains open. In the case that $H$ is not bipartite the exponent is 2, whereas
in the bipartite case the exponent is generally not known. Solymosi's conjecture is as follows:

\begin{conj}\label{soly}
If $H$ is a graph with exponent $\alpha$, then any $n$-vertex graph in which every edge is in exactly one copy of $H$ has $o(n^{\alpha})$ edges.
\end{conj}

In the case $H = C_4$ i.e.\ $H$ is a quadrilateral, Solym\'{o}si~\cite{Sol} conjectured that if an $n$ vertex graph is a union of $\Theta(n^{3/2})$ edge-disjoint quadrilaterals, then the graph contains $\Omega(n^2)$ quadrilaterals. By the removal lemma, Conjecture \ref{soly} is true for non-bipartite graphs $H$, so the conjecture is interesting for bipartite graphs.

\subsection{Main Result}

Our main result shows that there are infinitely many bipartite graphs $H$ for which Conjecture \ref{soly} is false.
Let $H_k$ be the graph with vertex set $V(H_k) = \{1,2,\dots,2k\}$ and edge set
$E(H_k) = \{1i,2i,3j,4j  : 3 \leq i \leq k + 2 < j \leq 2k\}$ -- see Figure~1.

\vspace{-0.3in}

\centerline{\includegraphics[scale=.4]{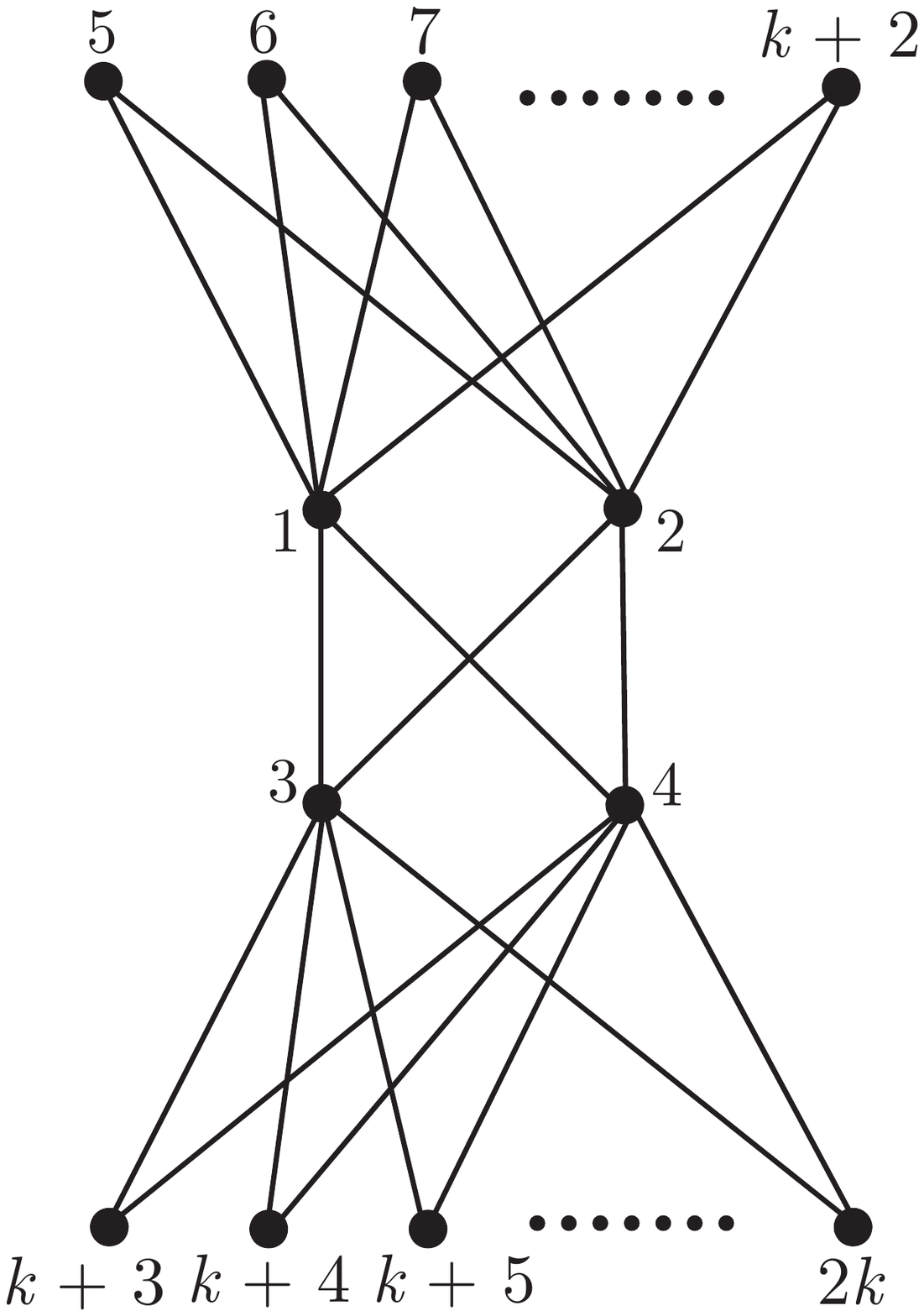}}

\vspace{-1.1in}

\begin{center}
Figure~1: The bipartite graph $H_k$.
\end{center}

We shall prove that $H_k$ has exponent $\frac{3}{2}$ for all $k \geq 3$ (Section \ref{exf}), and then the following theorem shows that the graphs $H_k$ for $k \geq 5$
give counterexamples to Conjecture \ref{soly}.

\begin{theorem}\label{main}
Let $k \geq 5$ and let $P$ be the set of primes $p \equiv 1$ mod 4. There exists a sequence of graphs
$(G_p)_{p \in P}$ such that $G_p$ has $n = 2kp^2$ vertices, $\Theta(n^{3/2})$ edges, and every edge of every $G_p$ is
contained in exactly one copy of $H_k$.
\end{theorem}

The proof of Theorem \ref{main} is based on a construction described in Section \ref{se:construction} involving simple field arithmetic and the quadratic character of $\fp$.  The proof is given
in Sections \ref{se:construction} -- \ref{exf}. The case that $H$ is a quadrilateral in Conjecture \ref{soly} remains open.

%%%%%%%%%%%%%%%%%%%%%%%%%%%%%%%%%%%%%%%%%%%%%%%%%%

\section{Outline of the Proof of Theorem 1.1}\label{se:construction}

The proof of Theorem \ref{main} is achieved in three steps. Before we describe these steps, we introduce some notation.
Throughout this section, $p$ is a prime, $k$ is a positive integer, and
$\fp$ denotes the finite field of order $p$. Let $\chi$ denote the quadratic character of $\fp$, namely $\chi(x) = 1$ if $x$ is a non-zero quadratic residue,
$\chi(x) = -1$ if $x$ is not a quadratic residue, and $\chi(0) = 0$.
The graph $H_k$ has vertex set $[2k] := \{1,2, \dots , 2k \}$ and edge set
\[ E(H_k) := \{1i,2i,3j,4j  : 3 \leq i \leq k + 2 < j \leq 2k\}.\]
Let $V(G)$ and $E(G)$ denote the vertex set and edge set of a graph $G$. If $G_1$ and $G_2$ are edge-disjoint graphs, 
then we write $G = G_1 \oplus G_2$ if $V(G) = V(G_1) \cup V(G_2)$ and $E(G) = E(G_1) \cup E(G_2)$.

\subsection{Step I. Construction of graphs $G_{\Gamma,\Lambda,S}(H)$.}

Let $\Gamma$ be a finite abelian group and $S \subseteq \Gamma$.
Let $H$ be an arbitrary graph with vertex set $[k]$, edge set $E$, and let $\Lambda = \{\lambda_1,\lambda_2,\dots,\lambda_k\} \subset \mathbb Z$.
For $i \in [k]$ let $X_i = \Gamma \times \{i\}$.  For $ij \in E$ with $i < j$, 
let $G_{\Gamma,\Lambda,S}(ij)$ be the bipartite graph with parts $X_i$ and
$X_j$, where $x \in X_i$ is adjacent to $y \in X_j$ whenever there exists $a \in S$
such that $y = x + (\lambda_j - \lambda_i)a$. Finally, define the following $k$-partite graph with parts $X_1,X_2,\dots,X_k$:
\[ G_{\Gamma,\Lambda,S}(H) = \bigcup_{ij \in E} G_{\Gamma,\Lambda,S}(ij).\]
  A key observation is that $G_{\Gamma,\Lambda,S}(H)$ is built up from edge-disjoint copies of $H$. The following is proved in Section \ref{basic}:

\begin{proposition}\label{union}
Let $k \in \mathbb N$ and let $H$ be a $k$-vertex graph.
For any finite abelian group $\Gamma$, $S \subseteq \Gamma$, and $\Lambda \subset \mathbb Z$, the graph $G_{\Gamma,\Lambda,S}(H)$ has $k|\Gamma|$ vertices, $|E||\Gamma||S|$ edges,
and there exist induced subgraphs $H(1),H(2),\dots,H(|\Gamma||S|)$ of $G_{\Gamma,\Lambda,S}(H)$, each isomorphic to $H$, such that
\[ G_{\Gamma,\Lambda,S}(H) = \bigoplus_{i = 1}^{|\Gamma||S|} H(i).\]
\end{proposition}

The proof of Theorem \ref{main} involves a suitable choice of the group $\Gamma$, the set $S$, and the set $\Lambda \subset \mathbb{Z}$. The graphs $G_{\Gamma,\Lambda,S}(H)$ 
may be useful as constructions for other extremal graph theory problems.

\subsection{Step II. The choice of $\Gamma$, $S$, and $\Lambda$.}

Let $\Gamma = \fp^2$ and $S = \{(a,a^2) : a \in \fp^*\}$ where $p \equiv 1$ mod 4 is prime. We aim to show that there exists a choice of $\Lambda = \{ \lambda_1 , \dots , \lambda_k \}$ such that every edge of
$G_p = G_{\Gamma,\Lambda,S}(H_k)$ is in exactly one copy of $H_k$, which is the heart of Theorem \ref{main}. It is convenient to let
$\lambda_{ij} = \lambda_j - \lambda_i$. After some arithmetic preparations in Section \ref{quad}, we prove the following
in Section \ref{suitproof}.

\begin{proposition}\label{suit}
Let $p \equiv 1$ mod 4 be a prime, $\Gamma = \fp^2$, and $S = \{(a,a^2) : a \in \fp^*\}$. Suppose $\Lambda = \{\lambda_i : i \in [2k]\}$ and
\begin{center}
\begin{tabular}{lp{5in}}
 {\rm 1.} & $\lambda_1 = 0$, $\lambda_2 = 1$. \\
 {\rm 2.} & For $i \in \{3,4 \}$ and $j \in \{1,2 \}$, $\chi( \lambda_{1i} \lambda_{i2} ) = \chi ( \lambda_{3j} \lambda_{j4} ) = -1$. \\
 {\rm 3.} & For $5 \leq i \leq k + 2 < j \leq 2k$, we have $\chi ( \lambda_{1i} \lambda_{i2} ) = \chi ( \lambda_{3j} \lambda_{j4} ) = 1$.
 \end{tabular}
 \end{center}
Then every edge of $G_p = G_{\Gamma,\Lambda,S}(H_k)$ is contained in exactly one copy of $H_k$.
\end{proposition}

Note that $\chi(\lambda_{ij}) = \chi(-\lambda_{ji}) = \chi(\lambda_{ji})$ since $p \equiv 1$ mod 4.
A set $\Lambda$ satisfying the conditions in Proposition \ref{suit} is called {\em suitable}.
We prove the following in Section \ref{choice}:

\begin{proposition}\label{lambdachoice}
If $k \geq 5$ and $p \geq 4k + 3$, then there exist a suitable set $\Lambda \subset \fp$.
\end{proposition}

If $\Gamma = \fp^2$ and $S = \{(a,a^2) : a \in \fp^*\}$, then it is not hard to show $G_{\Gamma,S,\Lambda}(C_4)$ contains
$\Theta(n^2)$ quadrilaterals for any set $\Lambda$, so $G_{\Gamma,S,\Lambda}(C_4)$ cannot be used as a counterexample to Conjecture \ref{soly} when $H = C_4$.

\subsection{Step III. The Tur\'{a}n Number for $H_k$.}

It remains to show $\ex(n,H_k) = \Theta(n^{3/2})$, which is the final step in the proof of Theorem \ref{main}.
A counting argument will be used to determine the order of magnitude of $\ex(n,H_k)$.

\begin{proposition}\label{exf}
For any integer $k \geq 3$, $\ex(n,H_k) = \Theta(n^{3/2})$.
\end{proposition}

This proposition is proved in Section \ref{exfsection}. Since $K_{2,k} \subset H_k$, the results of F\"{u}redi~\cite{Furedi} show
\[ \ex(n,H_k) \geq \ex(n,K_{2,k}) \geq \frac{1}{2} \sqrt{k-1}n^{3/2} + o(n^{3/2})\]
and this establishes the lower bound in Proposition \ref{exf}.

\section{Proof of Proposition 1}\label{basic}

It follows from the definition of $G_{\Gamma,\Lambda,S}(H)$ that $G_{\Gamma,\Lambda,S}(H)$ has $k|\Gamma|$ vertices.
We now prove that $G_{\Gamma,\Lambda,S}(H)$ is a union of edge-disjoint copies of $H$.
We observe that $H$ appears naturally as a subgraph of $G_{\Gamma,\Lambda,S}(H)$ in the following manner: for
$v \in \Gamma$ and $a \in S$, the set $\{v_i = v + \lambda_i a : i \in [k]\}$ induces a subgraph $H(v,a)$ isomorphic to
$H$, since $v_i v_j \in G_{\Gamma,\Lambda,S}(H)$ if and only if $ij \in E(H)$. Also, if $(x,i),(y,j) \in H(v,a)$ and $(x,i),(y,j) \in H(w,b)$, then
\[ x = v + \lambda_i a = w + \lambda_i b \quad y = v + \lambda_j a = w + \lambda_j b\]
and therefore $\lambda_i(a - b) = \lambda_j(a - b)$ which means $a = b$ and so $v = w$.
Therefore, no two subgraphs $H(v,a)$ share any pair of vertices. We conclude
\[ G_{\Gamma,\Lambda,S}(H) = \bigoplus_{(v,a) \in \Gamma \times S} H(v,a).\]
In particular, $G_{\Gamma,\Lambda,S}(H)$ has $|E||\Gamma||S|$ edges. This proves Proposition \ref{union}. \qed

\section{Proof of Proposition 2}\label{suitable}

This section is split into two parts.  First we describe the interaction between quadrilaterals in $G_p(H) := G_{\Gamma,\Lambda,S}(H)$ for general $H$, $\Gamma = \fp^2$, $\Lambda \subset \fp$ and $S = \{(a,a^2) : a \in \fp^*\}$. We then show that if the conditions of Proposition \ref{suit} are satisfied, then every edge of $G_p = G_p(H_k)$ is in exactly one copy of $H_k$.

\subsection{Quadrilaterals in $G_{\Gamma,\Lambda,S}(H)$.}\label{quad}

Throughout this subsection, $H$ is a graph, $G_p(H) = G_{\Gamma,\Lambda,S}(H)$ where $\Gamma = \fp^2$, $\Lambda \subset \fp$ and $S = \{(a,a^2) : a \in \fp^*\}$. 
The following simple but key arithmetic lemma is due to Ruzsa \cite{R}.

\begin{lemma}\label{Ruzsa}
Let $\alpha , \beta , \gamma , \delta,a,b,c,d \in \fp^*$ where $\alpha + \beta = \gamma + \delta \neq 0$. If
\begin{equation*}
\alpha (a,a^2) + \beta (b,b^2) = \gamma (c,c^2) + \delta (d , d^2),
\end{equation*}
then $\alpha \beta ( a - b)^2 = \gamma \delta (c - d)^2$.
\end{lemma}

\begin{pf}
Multiply $\alpha a^2 + \beta b^2 = \gamma c^2 + \delta d^2$ by $\alpha + \beta = \gamma + \delta$ to get
\begin{equation*}
\alpha^2 a^2 + \alpha \beta ( a^2 + b^2) + \beta^2 b^2 = \gamma^2 c^2 +
\gamma \delta (c^2 + d^2) + \delta^2 d^2.
\end{equation*}
Subtracting the square of $\alpha a + \beta b = \gamma c+ \delta d$ gives $\alpha \beta ( a-b)^2 = \gamma \delta ( c-d)^2$.
\end{pf}

This lemma has a number of consequences relative to the distribution of quadrilaterals in $G(ij) := G_{\Gamma,\Lambda,S}(ij)$ for $ij \in H$.
For $hi,ij \in E(H)$, let $G(hij) = G(hi) \cup G(ij)$.

\begin{lemma}\label{c4}
For any edge $ij \in E(H)$, the graph $G(ij)$ is $C_4$-free.
\end{lemma}

\begin{pf}
Consider a quadrilateral $(x,i)(y,j)(z,i)(w,j) \subset G(ij)$.
For some $a,b,c,d \in \fp^*$,
\[ y = x + \lambda_{ij} (a,a^2) = z + \lambda_{ij} (b,b^2) \quad w = x + \lambda_{ij} (c,c^2) = z + \lambda_{ij} (d,d^2).\]
Canceling out $x,y,z$ and $w$ and dividing by $\lambda_{ij}$, we obtain in each component
$a + b = c + d$ and $a^2 + b^2 = c^2 + d^2$. By Lemma \ref{Ruzsa}, $(a - b)^2 = (c - d)^2$
and so $a - b = c - d$ or $a - b = d - c$. In either case, together with $a + b = c + d$ we find $\{a,b\} = \{c,d\}$.
But then $(x,i) = (z,i)$ or $(y,j) = (w,j)$, a contradiction. Therefore $G(ij)$ has no
$C_4$.
\end{pf}

\begin{lemma}\label{k23}
For any edges $hi,ij \in E(H)$, the graph $G(hij)$ is $K_{2,3}$-free.
\end{lemma}

\begin{pf}
Suppose $G(hij)$ contains a $K_{2,3}$. By Lemma \ref{c4}, $G(hi)$ and $G(ij)$ are $C_4$-free,
and so the three vertices of degree two in the $K_{2,3}$ must be in $X_i$, say
$(z_1,i),(z_2,i),(z_3,i)$, and the other two vertices are $(x,h) \in X_h$ and $(y,j) \in X_j$.
By definition, for some $a_1,a_2,a_3 \in \fp^*$ and $b_1,b_2,b_3 \in \fp^*$,
\[ y = x + \lambda_{hi}(a_r,a_r^2) + \lambda_{ij}(b_r,b_r^2)\]
for $r\in \{ 1,2,3 \}$.
Note that the $a_r$ are distinct and the $b_r$ are distinct, for if $a_r = a_s$ or $b_r = b_s$ for some $r \neq s$,
then $z_r = z_s$ and so $(z_r,i) = (z_s,i)$, a contradiction. On the other hand, by Lemma \ref{Ruzsa},
$(a_r - b_r)^2 = (a_s - b_s)^2$ for all $r,s \in \{1,2,3\}$. Taking square roots,
some pair of square roots has the same sign, namely for some $r \neq s$, we have $a_r - b_r = a_s - b_s$. For all $r,s$, we also have
\[ \lambda_{hi} a_r + \lambda_{ij} b_r = \lambda_{hi} a_s + \lambda_{ij} b_s.\]
Subtracting $\lambda_{hi}(a_r - b_r) = \lambda_{hi}(a_s - b_s)$ from this equation,
we obtain $\lambda_{hj} b_r = \lambda_{hj} b_s$.  Therefore, $b_r = b_s$ and $(z_r , i )  = (z_s , i)$ which is a contradiction. 
\end{pf}

\begin{lemma}\label{absorb1}
If $C = (x,i)(y,j)(z,k)(w,l)$ is a quadrilateral in $G_p(H)$
and $(x,i)(y,j)(z,k) \subset H(v,a)$, then $C \subset H(v,a)$.
\end{lemma}

\begin{pf}
By definition of $G_p(H)$, there exist $a,c,d \in \fp^*$ such that
\[ y = x + \lambda_{ij}(a,a^2) \quad z = y + \lambda_{jk}(a,a^2) \quad w = z + \lambda_{kl}(c,c^2) \quad x = w + \lambda_{li}(d,d^2).\]
This implies
\[ \lambda_{ij}(a,a^2) + \lambda_{jk}(a,a^2) = \lambda_{lk}(c,c^2) + \lambda_{il}(d,d^2).\]
By Lemma \ref{Ruzsa} with $\alpha = \lambda_{ij}$, $\beta = \lambda_{jk}$, $\gamma = \lambda_{lk}$, and
$\delta = \lambda_{il}$, noting $\alpha + \beta = \gamma + \delta$, we have
\[ \gamma\delta (c - d)^2 = 0\]
and therefore $c = d$. Now
\[ \lambda_{ik} a = \lambda_{ij}a + \lambda_{jk}a = \lambda_{lk}c + \lambda_{il} d = \lambda_{ik} c\]
so we conclude $a = c = d$. In particular, since $z = v + \lambda_{0k}(a,a^2)$,
\[ w = z + \lambda_{kl}(a,a^2) = v + \lambda_{0k}(a,a^2) + \lambda_{kl}(a,a^2) = v + \lambda_{0l}(a,a^2)\]
and we conclude $C \subset H(v,a)$.
\end{pf}

\begin{lemma}\label{trivialc4}
If $C = (x,i)(y,j)(z,k)(w,l)$ is a quadrilateral in $G_p(H)$ and
$\chi(\lambda_{ij}\lambda_{jk}\lambda_{kl}\lambda_{li}) = -1$ then $C \subset H(v,a)$ for some $(v,a) \in \fp^2 \times S$.
\end{lemma}

\begin{pf}
By definition of $G_p(H)$, there exist $a,b,c,d \in \fp^*$ such that
\[
\lambda_{ij} (a,a^2) + \lambda_{jk} (c,c^2) = \lambda_{il} (d,d^2) + \lambda_{lk} (b,b^2).
\]
By Lemma 1,
\[
\lambda_{ij} \lambda_{jk} (a-c)^2 = \lambda_{il} \lambda_{lk} (d - b)^2.
\]
Since $\chi( \lambda_{ij} \lambda_{jk} \lambda_{il} \lambda_{lk} ) 
= \chi ( \lambda_{ij} \lambda_{jk} \lambda_{kl} \lambda_{li} ) = -1$, we conclude 
$a = c$ and $b = d$.  The equation 
$\lambda_{ij}a + \lambda_{jk}c  = \lambda_{il}d + \lambda_{lk}b$ 
reduces to $\lambda_{ik}a = \lambda_{ik}d$ and we conclude $a=b=c=d$.
Letting $v = x - \lambda_{0i} (a,a^2)$, we have $C \subset H(v,a)$.  
\end{pf}

\subsection{Proof of Proposition 2}\label{suitproof}

Suppose $F \subset G_p = G_{ \Gamma , \Lambda , S} (H_k)$ is isomorphic to $H_k$ and let $\phi : V(H_k) \rightarrow V(F)$ be an isomorphism. We aim to show that $F = H (v,a)$ for some $v \in \fp$ and $a \in S$ by finding, via Lemma \ref{trivialc4}, a quadrilateral
$C^* \subset H_k$ with $\phi(C^*) \subset H(v,a)$. This is the point where we make heavy use
of the last two conditions in Proposition \ref{suit}. Let $\mathfrak{c}: V(H_k) \rightarrow \Lambda$ be the proper vertex-coloring of $H_k$ given by
$\mathfrak{c}(x) = \lambda_i$ if and only if $\phi(x) \in X_i$. Let $[m,n] := \{m,m+1,m+2,\dots,n\}$. 

\medskip

{\bf Claim 1.} {\em Each quadrilateral in $H_k$ is colored with at least three colors, and each $K_{2,3}$ in $H_k$
is colored with at least four colors. Furthermore, }
\begin{equation*}\label{color eq}
\{ \mathfrak{c}(1) , \mathfrak{c}(2) , \mathfrak{c}(3) , \mathfrak{c}(4) \} \subset \{\lambda_1 , \lambda_2 , \lambda_3, \lambda_4 \}.
\end{equation*}

{\it Proof of Claim 1.} The first statement is an immediate consequence of Lemmas \ref{c4} and \ref{k23}.
Consider the vertex $1 \in V(H_k)$ (see Figure~1).
If $\mathfrak{c}(1) \in \{ \lambda_5 , \lambda_6 , \dots , \lambda_{k+2}\}$
then each neighbor of $1$ is assigned color $\lambda_1$ or color $\lambda_2$.  
This follows from the fact that if $xy \in E(H_k)$ and $\mathfrak{c}(x) = \lambda_i$ and $\mathfrak{c}(y) = \lambda_j$, then 
$ij \in E(H_k)$.  
The neighbors of 1 are also neighbors of 2 so that $\mathfrak{c}$ assigns one color to at least three common neighbors of 1 and 2. This is impossible by the first part of the claim.  A similar argument shows that
$\mathfrak{c}(1) \notin \{ \lambda_{k+3} , \lambda_{k + 4} , \dots , \lambda_{2k} \}$ thus
$\mathfrak{c}(1) \notin \{ \lambda_5 , \lambda_6 , \dots , \lambda_{2k}\}$.
By symmetry, we must have
$\mathfrak{c}(i) \notin \{\lambda_5, \lambda_6,  \dots , \lambda_{2k}\}$ for $1 \leq i \leq 4$ and so
\begin{equation*}
\{ \mathfrak{c}(1) , \mathfrak{c}(2) , \mathfrak{c}(3) , \mathfrak{c}(4) \} \subset \{\lambda_1 , \lambda_2 , \lambda_3, \lambda_4 \}.
\end{equation*}
This proves Claim 1. \qed

\bigskip

{\bf Claim 2.} {\em There is a quadrilateral $C^* = ghij$ in $H_k$ such that}
\[
\chi \left( (  \mathfrak{c}(g) - \mathfrak{c}(h)  )( \mathfrak{c}(h) - \mathfrak{c}(i)   ) ( \mathfrak{c}(i) - \mathfrak{c}(j)  ) ( \mathfrak{c}(j) - \mathfrak{c}(g)  )  \right) = -1.
\]
{\it Proof of Claim 2.} By Claim 1 there is a path $ghi \subset 1324$ such that $\mathfrak{c}(g),\mathfrak{c}(h),\mathfrak{c}(i)$ are distinct. 
Without loss of generality, assume that $g=1$, $h=3$, and $i =2$.  
If we can find a $j \in [5,k+2]$ such that $\mathfrak{c}(j) \not \in \{\lambda_1,\lambda_2,\lambda_3,\lambda_4\}$ then by
the conditions of Proposition \ref{suit},
\[ \chi(\lambda_{13}\lambda_{32}\lambda_{2j}\lambda_{j1}) = \chi(\lambda_{13}\lambda_{32}) \cdot \chi(\lambda_{2j}\lambda_{j1}) = -1\]
and so $132j$ is the required cycle. To find $j$, note that no two vertices in $[5,k+2]$ 
have color $\mathfrak{c}(3)$  
otherwise those two vertices together with $1,3$ and $2$ form a 3-colored $K_{2,3}$, contradicting Claim 1. 
Similarly, no two vertices in $[5,k+2]$ have color $\mathfrak{c}(4)$.
Since $k \geq 5$, there are at least three vertices in $[5,k+2]$ and so one of them, say $j$, has a color not in $\{\lambda_1,\lambda_2,\lambda_3,\lambda_4\}$, as required. \qed

\medskip

We now complete the proof of Proposition \ref{suit}.
By Lemma \ref{trivialc4}, there exists $v \in \fp^2$ and $a \in S$ such that $\phi(C^*) \subset H(v,a)$. Given any edge $e$ of $F = \phi(H_k)$,
there exist quadrilaterals $C^* = C^1,C^2,\dots,C^r \subset H_k$ such that $e \in E(\phi(C^r))$, and $E(\phi(C^i)) \cap E(\phi(C^{i + 1})) \neq \emptyset$ for $i < r$.
By Lemma \ref{absorb1}, we inductively have $\phi(C^i) \subset H(v,a)$ for all $i \leq r$. In particular, $e \in E(H(v,a))$ and we conclude 
$F \subset H(v,a)$. Since $H(v,a)$ is an induced subgraph of $G_p(H_k)$ isomorphic to $F$, we conclude $F = H(v,a)$. \qed

\section{Proof of Proposition 3}\label{choice}

To build a suitable set $\Lambda \subset \fp$, we use the following identity (see Theorem 5.48, \cite{LN}):

\begin{proposition}\label{char} Let $f(x) = a_2 x^2 + a_1 x + a_0 \in \mathbb{F}_p [x]$ where $a_2 \neq 0$.
If $a_1^2 - 4a_0 a_2 \neq 0$ then
\[
\sum_{c \in \mathbb{F}_p} \chi (f(c)) = - \chi(a_2).
\]
\end{proposition}

Let $k \geq 5$, $p \geq 4k + 3$ be prime with $p \equiv 1 (\textup{mod}~4)$, and $\lambda_1 = 0$, $\lambda_2 =1$. By Proposition \ref{char},
\begin{equation}\label{tech lemma eq1}
\sum_{c \in \mathbb{F}_p} \chi ( c^2 - c) = -1
\end{equation}
so there are least $\frac{p-3}{2} \geq 2k$ elements of $\mathbb{F}_p$ for which $\chi (c^2 - c) = -1$.  Let $\lambda_3$ be any one of them and observe that since $\chi( \lambda_3^2 - \lambda_3 ) = -1$, we have $\lambda_3 \neq 0$ and $\lambda_3 \neq 1$.  Using the fact that
$\chi (-1) = 1$,
\[
-1 = \chi ( \lambda_3^2 - \lambda_3 ) = \chi ( (0 - \lambda_3) ( \lambda_3 - 1) ) = \chi ( (\lambda_1 - \lambda_3) ( \lambda_3 - \lambda_2) ).
\]

Next we choose $\lambda_4$.  Let $g(x) = (1 - \chi (x^2 - x) )( 1 - \chi ( \lambda_3 x) )$ and
$X = \{ c \in \mathbb{F}_p : g(c) = 4 \}$.  If $c \in X$ then $\chi (c^2 - c) = \chi (\lambda_3 c) = -1$.  Suppose $c \in \mathbb{F}_p$ and $0 < g(c) \leq 2$.  Then either $\chi (c^2 - c) = 0$ or
$\chi ( \lambda_3 c ) = 0$ which means $c = 0$ or $c = 1$ thus
\begin{equation}\label{tech lemma eq2}
\Bigl| \sum_{c \in \mathbb{F}_p} g(c)  - 4 |X| \Bigr| \leq 2 \cdot 2.
\end{equation}
Expanding $g(c)$ and using Proposition \ref{char},
\begin{equation}\label{tech lemma eq3}
\sum_{c \in \mathbb{F}_p} g(c) = p + \sum_{c \in \mathbb{F}_p} \left(
- \chi (c^2 - c) - \chi(\lambda_3) \chi (c) + \chi(c^2) \chi (\lambda_3) \chi (c - 1) \right) = p +1.
\end{equation}
Here we have used the well known fact that $\sum_{c \in \mathbb{F}_p} \chi (c) = 0$.
By (\ref{tech lemma eq2}) and (\ref{tech lemma eq3}),
\[
|X| \geq \frac{ p +1}{4} - 1
\]
and since $\frac{p-3}{4} \geq 1$, the set $X$ is not empty.  Choose $\lambda_4$ so that
$-1 = \chi ( \lambda_4 ( \lambda_4 - 1) ) = \chi (\lambda_3 \lambda_4)$.  Then $\lambda_4$ is not equal to any of $\lambda_1  ,\lambda_2$, or $\lambda_3$.
The relation $-1 = \chi ( \lambda_4 ( \lambda_4 -1) )$ implies $-1 = \chi ( ( \lambda_1  - \lambda_4)(\lambda_4 - \lambda_2))$ and
the relation $-1 = \chi (\lambda_3 \lambda_4)$ implies $-1 = \chi ( (\lambda_3 - \lambda_1)(\lambda_1 - \lambda_4))$.  Furthermore
\begin{eqnarray*}
\chi ( ( \lambda_3 - \lambda_2)(\lambda_2 - \lambda_4) ) &  = & \chi ( ( \lambda_3 -1)(\lambda_4 - 1) )\\
& =&  \chi ( \lambda_3 (\lambda_3 - 1) ) \chi (\lambda_3 \lambda_4) \chi ( \lambda_4 (\lambda_4 - 1) ) \\
& = & (-1)(-1)(-1) = -1
\end{eqnarray*}
so that all of the requirements of Condition 2 are satisfied.

Choosing the remaining $\lambda_i$'s will be straightforward.  By (\ref{tech lemma eq1}), there are at least
$\frac{p-3}{2}$ elements $c \in \mathbb{F}_p$ for which $\chi (c^2 -c) = 1$.  Since
$\frac{p-3}{2} \geq 2k$ we can choose $\lambda_5, \lambda_6 , \dots , \lambda_{k+2}$ so that none of these are equal to
$\lambda_1$, $\lambda_2$, $\lambda_3$ or $\lambda_4$, and $\chi( (  \lambda_1 - \lambda_i ) (\lambda_i - \lambda_2) )  = 1$ for all $4 < i \leq k+2$.

Let $f(x) = (\lambda_3 - x)(x-\lambda_4) = - x^2 + (\lambda_3 + \lambda_4)x - \lambda_3 \lambda_4$.  If $a_2 = -1$, $a_1 = \lambda_3 + \lambda_4$, and
$a_0 = -\lambda_3 \lambda_4$, then $a_1^2 - 4a_0 a_2 = ( \lambda_3  - \lambda_4)^2 \neq 0$.  By Proposition \ref{char},
\[
\sum_{c \in \mathbb{F}_p} \chi ( f(c) ) = - \chi ( -1) = -1
\]
and there are at least $\frac{p-3}{2}$ elements $c$ for which $\chi (f (c) ) = 1$.  Since
$\frac{ p -3}{2} \geq 2k$ we can choose $\lambda_{k+3} , \lambda_{k+4} , \dots , \lambda_{2k}$ so that
all of $\lambda_1 , \lambda_2 , \dots , \lambda_{2k}$ are distinct and for any $4 < i \leq k+2$,
\[
\chi ( ( \lambda_3 - \lambda_{i+k - 2} )( \lambda_{i+k - 2} - \lambda_4) ) = 1.
\]
This completes the proof of Proposition \ref{choice}. \qed

\section{Proof of Proposition 4}\label{exfsection}

The claim $\ex(n,H_k) = \Theta(n^{3/2})$ follows by proving $\ex(n,H_k) \leq kn^{3/2}$ for all $k \geq 3$ and large enough $n$.
We follow the method of dependent random choice~\cite{FoxSudakov}, but we do not optimize the constants in the upper bound we obtain for $\ex(n,H_k)$.

\medskip

Let $G$ be an $H_k$-free $n$-vertex graph with average degree $d := 2|E(G)|/n$. It is sufficient to show $d \leq 2k\sqrt{n}$.
Choose uniformly at random a pair of vertices $\{x_1,x_2\}$ in $G$ and let $S = N(x_1,x_2)$, the common neighborhood of
$x_1$ and $x_2$.  Given a pair of vertices $\{y_1 , y_2 \}$, let $d( y_1 ,y_2) = |N(y_1,y_2)|$.
Let $X = |S|$ and let $Y$ be the number of $\{y_1,y_2\} \subset S$ with $d(y_1,y_2) \leq 2k$.
If $X - Y \geq k + 2$, then there exist two vertices $\{x_1,x_2\} \subset V(G)$ with $d(x_1,x_2) > k+2$ and
some pair $\{y_1,y_2\} \subset N(x_1,x_2)$ with $d(y_1,y_2) > 2k$, and we easily find a copy of $H_k$ by mapping
$\{1,2\}$ to $\{x_1,x_2\}$ and $\{3,4\}$ to $\{y_1,y_2\}$, a contradiction. So $X - Y \leq k+1$.

\medskip

On the other hand, using convexity of binomial coefficients,
\begin{eqnarray*}
k + 1 &\geq& \mathbb E(X - Y) \\
&\geq& \sum_{y \in V(G)} \frac{\binom{d(y)}{2}}{\binom{n}{2}} -
\binom{n}{2} \cdot \frac{\binom{2k}{2}}{ \binom{n}{2} } \\
& \geq&  \frac{1}{ \binom{n}{2} } n \binom{d}{2} - \binom{2k}{2} \\
&\geq& \frac{d(d-1)}{n}  - \binom{2k}{2}.
\end{eqnarray*}
It follows that if $n$ is large enough, then $d \leq 2k\sqrt{n}$. This proves Proposition \ref{exf}.
\qed

\section{Concluding remarks}

$\bullet$ A {\em Sidon set} in a finite abelian group $\Gamma$ is a set $S \subset \Gamma$ such that if $a + b = c + d$ with $a,b,c,d \in S$, then $\{a,b\} = \{c,d\}$.
A generalization of Sidon sets was given in~\cite{LV}.  Consider all equations $\alpha a + \beta b = \gamma c + \delta d$
where $\alpha + \beta = \gamma + \delta$ and $1 \leq \alpha,\beta,|\gamma|,|\delta| \leq k$. If $A \subset \Gamma$ has no solution
to any of these equations other than $\{a,b\} = \{c,d\}$, then $A$ is called a {\em $k$-fold Sidon set}. It is easy to see that a
$k$-fold Sidon set in a finite abelian group $\Gamma$ has size $O(|\Gamma|^{1/2})$ (see~\cite{CT} for more details).
A 2-fold Sidon set of size roughly $\sqrt{N}/2$ in $\mathbb Z_N$ is constructed for infinitely many $N$ in~\cite{LV}, but it is an open question to construct a $k$-fold Sidon set
of size $\Theta(\sqrt{N})$ in any abelian group $\Gamma_N$ of order $N$ for any $k \geq 3$. In fact the following is conjectured in~\cite{LV}:

\begin{conj}
Let $k \in \mathbb N$. Then there exists $c_k > 0$ such that for any $N \in \mathbb N$, there exists a 
$k$-fold Sidon set $A \subset [N]$ of size at least $c_k\sqrt{N}$. 
\end{conj}

The densest current construction available is due to Ruzsa~\cite{R}, who showed that for each $k$ there exists a $k$-fold Sidon
set of size $N^{1/2 - o(1)}$ in $\mathbb Z_N$. The construction used for Theorem \ref{main} is easily adapted to give counterexamples
to Conjecture \ref{soly} if there exists a $k$-fold Sidon set of size $\Theta(\sqrt{N})$ in an abelian group $\Gamma_N$ of order $N$.
The following theorem is proved in the same way as Theorem \ref{main}, by taking $G_N = G_{\Gamma_N,\Lambda,S}(C)$ where $C$ is a quadrilateral,
$\Lambda = \{0,3,1,2\}$ and $S$ is a 3-fold Sidon set:

\begin{theorem}\label{kfoldthm}
If there exists a $3$-fold Sidon set of size $\Theta(\sqrt{N})$ in an abelian group $\Gamma$ of order $N$,
then there exists a $4N$-vertex graph $G_N$ with $\Theta(N^{3/2})$ edges such that every edge of $G_N$ is contained in exactly one quadrilateral.
\end{theorem}

A particularly case is to find a Sidon set of size $\Theta(\sqrt{N})$ in $\mathbb Z_N$ such that no difference of two distinct elements equals twice or three times the difference of any other two elements; i.e. the equations $a - b = 2(c - d)$ and $a - b = 3(c - d)$ force $a = b$ and $c = d$.

\medskip

$\bullet$ A more general setting is given in~\cite{LV} in the language of hypergraphs~\cite{Berge}.
An $r$-uniform hypergraph has {\em girth five} if whenever one pair of vertices is selected from each hyperedge, the resulting
graph has no cycles of length at most four (and in particular no double edges). In~\cite{LV} it is conjectured that for every $r \geq 2$, there exists an
$r$-uniform hypergraph on $n$ vertices with girth five and $\Theta(n^{3/2})$ hyperedges. This is settled for $r = 2$ by Erd\H{o}s and R\'{e}nyi~\cite{ER}, and for $r = 3$ in~\cite{LV}.
 In the case $r = 4$, if this conjecture is true then we may place
a quadrilateral in each hyperedge to obtain a graph with $\Theta(n^{3/2})$ edges and $n$ vertices, in which every edge is in exactly one quadrilateral.
Using Ruzsa's construction~\cite{R}, one finds for each $r \geq 3$ an $n$-vertex $r$-uniform hypergraph of girth five with $n^{3/2 - o(1)}$ edges.
In particular, if $f_r(n)$ is the maximum number of hyperedges in an $r$-uniform $n$-vertex hypergraph of girth five, then 
$f_2(n) = \Theta (n^{3/2})$, $f_3(n) = \frac{1}{6} n^{3/2}+ o(n^{3/2})$, and 
for some constant $c_r > 0$,
\[ \frac{n^{3/2}}{\exp(c_r\sqrt{\log n})} \leq f_r(n) \leq \frac{1}{r(r - 1)}n^{3/2} + O(n)\]
for all $r \geq 4$.

%$\bullet$ In fact, if $f(x)$ is any polynomial with integer coefficients, it is possible to show that the equations
%\[  \alpha a + \beta b = \gamma c + \delta d \quad \tilde{\alpha} f(a) + \tilde{\beta} f(b) = \tilde{\gamma} f(c) + \tilde{\delta} f(d)\]
%where $\alpha + \beta = \gamma + \delta$, $\alpha \neq -\beta$, and $\tilde{\alpha} + \tilde{\beta} = \tilde{\gamma} + \tilde{\delta}$ has
%$(a,b,c,d)$ with $\{a,b\} \neq \{c,d\}$, and so it seems that the use of a set $S = \{(a,f(a) : a \in \fp^*\}$ cannot be used to construct an appropriate
%graph where each edge is in exactly one quadrilateral.

\end{document}